%% file: colouredCFG.tex
\documentclass{article}

\usepackage{theorem}
\usepackage{latexsym}
\usepackage{graphicx}
\graphicspath{{figures/}}
\usepackage{gastex}

{
\theorembodyfont{\rmfamily}
\newtheorem{example}{Example}[section]
\newtheorem{note}[example]{Note}
}

\newtheorem{theoreme}[example]{Theorem}
\newtheorem{corollaire}[example]{Corollary}

\newtheorem{definition}[example]{Definition}
\newtheorem{proposition}[example]{Proposition}

\newtheorem{lemme}[example]{Lemma}

\newcommand{\qed}{\hfill$\Box$ \vspace{0.5 cm}}

\newenvironment{proof}{{\it Proof~: }}{\qed \vskip 5mm}

\newcommand{\mylabel}[1]{\label{#1}}

\newcommand{\lcfg}{\ensuremath{L(CFG)}}

\newcommand{\lc}{\ensuremath{L(C)}}
\newcommand{\confspace}[1]{\ensuremath{L(#1)}}
\newcommand{\shot}{shot-set}
\newcommand{\sh}[1]{\ensuremath{s(#1)}}
\newcommand{\J}[1]{\ensuremath{J_{#1}}}
\newcommand{\M}[1]{\ensuremath{M_{#1}}}
\newcommand{\liste}[2][n]{\ensuremath{#2_1,\ldots,#2_{#1}}}
\newcommand{\ie}{\emph{i.e.}}
\def\N{\mbox{I\hspace{-.15em}N}}

\title{Characterisation of Lattices Induced by (extended) Chip Firing Games}
\author{Cl\'emence Magnien \and Ha Duong Phan \and Laurent Vuillon
	\thanks{(magnien,phan,vuillon)@liafa.jussieu.fr}}
\date{}

\begin{document}
\maketitle

\noindent
\textbf{Abstract:}
The Chip Firing Game (CFG) is a discrete dynamical model used in physics, computer science and economics.
It is known that the set of configurations  reachable from an initial configuration (this set is called the \emph{configuration space}) can be ordered as a lattice.
We first present a structural result about this model,
which allows us to introduce some useful tools for describing those lattices.
Then
we establish that the class of lattices that are the configuration space of a CFG is strictly between the class of distributive lattices and the class of upper locally distributive (or ULD) lattices.
Finally we propose an extension of the model, the \emph{coloured} Chip Firing Game, which generates exactly the class of ULD lattices.\\

\noindent
\textbf{keywords:}
Chip Firing Game, lattice, discrete dynamical model, Sand Pile Model

\section{Introduction}

The Chip Firing Game (CFG) was introduced by Bj\"orner, Lovasz and Shor in \cite{BLS91} and \cite{BL92}.
It is defined over a directed multigraph $G=(V,E)$, called the \emph{support graph} of the game.
A \emph{configuration} of the game is 
a mapping $\sigma:V\mapsto \N$ that associates
 a weight to each vertex, which can be considered as a number of \emph{chips} stored in the vertex.
The CFG is a discrete dynamical model, with the following evolution rule, also called the \emph{firing} rule:
if, when the game is in a configuration $\sigma$, a vertex $v$ contains at least as many chips as its outgoing degree,
one can transfer a chip from $v$ along each of its outgoing edges to the corresponding vertex.
We also call this \emph{applying the rule $v$} (and we will speak equivalently of $v$ as a vertex or as a \emph{rule}). If $\sigma'$ is the resulting 
configuration, we denote it by: $\sigma\stackrel{v}{\longrightarrow} \sigma'$, and we call $\sigma$ a \emph{predecessor} of $\sigma'$.

CFGs are strongly convergent games \cite{Eri93}, which means that, given an initial configuration, 
either the game can be played forever,
or it reaches a unique fixed point (where no firing is possible)
 independent on the order in which the vertices were fired. 
We will consider here only CFGs that reach a fixed point, that we call the
 \emph{final configuration} of the CFG.
These CFGs are \emph{convergent} CFGs.
(it is possible to guarantee that a CFG is convergent by the presence in the support graph of a sink accessible from all vertices \cite{LP00}).
We call \emph{execution} of a CFG any sequence of firing that, from the initial configuration, reaches the final configuration.

We know that
the configuration space of a convergent CFG, ordered by the reflexive and transitive closure of the predecessor relation, is a lattice \cite{LP00} (for an introduction to lattice theory see \cite{DP90}).
Moreover, this lattice is ranked, which means that all the paths from the initial configuration to the final configuration have the same length.

Given a CFG $C$, we denote by \confspace{C} its configuration space considered as a lattice.
We denote by \emph{\lcfg{}} the class of lattices that are the configuration space of a CFG.
Given a lattice $L\in\lcfg$,
if a CFG  $C$ is such that $\confspace{C}=L$, we say that $C$ is a 
CFG \emph{corresponding} to $L$.
We say that two CFGs are \emph{equivalent} if the lattices of their configuration spaces are isomorphic.

In Section \ref{secrecalls}, we present the definitions and results that are needed to study the lattices in \lcfg{}.
In Section \ref{seccfgsimple} we show that any CFG is equivalent to a simple CFG, \ie{} a CFG where each vertex is fired at most once during an execution,
 and we introduce efficient tools to describe the configuration spaces of such CFGs.
In Section \ref{secstudy} we attempt to characterise the class \lcfg{}.
In Section \ref{seccoloured} we present an extension of the model, the coloured Chip Firing Game.

\section{Recalls and definitions} \mylabel{secrecalls}

Let us first recall some basic definitions about posets and lattices \cite{DP90}:
a partially ordered set (or poset) is a set equipped with an order relation $\le$ (\ie{} transitive, reflexive and antisymmetric).
A \emph{linear extension} of a poset $P$ is a list \liste{x} of all its elements such that $i<j$ implies $x_i<x_j$.
An \emph{ideal} of a poset $P$ is a subset $I$ of $P$ such that, for all $x\in I, y\le x$ implies $y\in I$.
We will denote by $\cal{O}(P)$ the set of all ideals of $P$.
If $x$ and $y$ are two elements of a poset, we say that $x$ is \emph{covered} by $y$ (or $y$ \emph{covers} $x$), and write $x\prec y$ or $y\succ x$, 
if $x<y$ and $x\le z<y$ implies $z=x$.
To represent a poset $P$ we will use its Hasse diagram, defined as follows~:
\begin{itemize}
\item each element $x$ of $P$ is represented by a point $p_x$ of the plane,
\item if $x<y$, then $p_x$ is lower than $p_y$,
\item $p_x$ and $p_y$ are joined by a line if and only if $x\prec y$.
\end{itemize}

A poset $L$ is a \emph{lattice} if any two elements $x$ and $y$ of $L$ have a 
least upper bound (called \emph{join} and denoted by $x\vee_L y$ or simply $x\vee y$)
 and a greatest lower bound (called \emph{meet} and denoted by $x\wedge_L y$ or $x\wedge y$).
All the lattices considered in this paper are finite, therefore they have a least and a greatest element, respectively denoted by $0_L$ and $1_L$.
A subset $L_1$ of a lattice $L$ is a sub-lattice of $L$ if $L_1$ is stable by the join and meet in $L$: that is if, for each $x,y\in L_1, x\vee_L y\in L_1$ and $x\wedge_L y\in L_1$.
A lattice is \emph{ranked} if all the paths in the covering relation from the maximal to the minimal element have the same length.

A lattice $L$ is \emph{distributive} if it satisfies one of the two following laws of distributivity (which are equivalent):
$$\forall x,y,z \in L,\  x\wedge(y\vee z)=(x\wedge y)\vee(x\wedge z)$$
$$  \forall x,y,z \in L,\  x\vee(y\wedge z)=(x\vee y)\wedge(x\vee z)$$

A lattice is a \emph{hypercube of dimension $n$} if it is isomorphic to the set of all subsets of a set of $n$ elements, ordered by inclusion.
It is also called a \emph{boolean lattice}.

A lattice is \emph{upper locally distributive} (denoted by \emph{ULD} \cite{Mon90}) if the interval between an element and the meet of all its upper covers is a hypercube.

An element $j$ of a lattice $L$ is a \emph{join-irreductible} of $L$ if it is not the join of any subset of $L$ that does not contain $j$. 
Dually, $m\in L$ is a \emph{meet-irreductible} if it is not the meet of any subset of $L$ that does not contain $m$. 
The join-irreductibles and meet-irreductibles of a lattice are easily recognisable in the diagram of a lattice,
since we have the following characterisation:
\begin{itemize}
\item $j$ is a join-irreductible if and only if it has a unique lower cover, denoted by $j^-$.
\item $m$ is a meet-irreductible if and only if it has a unique upper cover, denoted by $m^+$.
\end{itemize}
The set of join-irreductibles of a lattice $L$ is denoted by $J_L$, or simply $J$.
The set of meet-irreductibles is denoted by $M_L$, or $M$.

Let us recall here some definitions and results about lattices (unless explicitely specified, they come from \cite{Cas98}):
\begin{proposition}
Let $L$ be a lattice. Any element $x$ of $L$ is the join of the join-irreductibles that are smaller than itself, and the meet of the meet-irreductibles that are greater than itself:
$$x= \bigvee \{j \in J,j \leq x\}= \bigwedge \{m \in M, x \leq m\}$$
\end{proposition}
We denote by $J_x$ (resp. $M_x$) the set $\{j \in J$: $j \leq x\}$
(resp. $\{m \in M$: $x \leq m\}$). 
These sets are a coding of the lattice \cite{BF48,BM70}.
Indeed, for any elements $x,y$ of a lattice, the order relation is characterised by:
$$ x\leq y\iff J_x\subset J_y \iff M_y\subset M_x$$
Moreover, in a lattice, the join is given by the following formula \cite{BF48,BM70}:
$$M_{x\vee y}=M_x\cap M_y$$
For ULD lattices, we have the following characterisation:
\begin{proposition}
A lattice $L$ is ULD if and only if, for all $x,y\in L$,
$$y\succ x \iff |M_x\backslash M_y|=1$$
\end{proposition}
As a consequence, an ULD lattice is ranked, and its height is equal to $|M|$.
This also allows us to associate to each edge $(x,y)$ of an ULD lattice the meet-irreductible $m$ such that
 $\{m\}=M_x\backslash M_y$. We label the edge $(x,y)$ by $m$.

We now introduce the \emph{arrow relations}, useful for the proofs and characterisations:
\begin{definition} \cite{WIL83}
Let $L$ be a lattice, $j\in J$ and $m\in M$. We define:
\begin{itemize}
\item $j\downarrow m$ if $j\not\le m$ and $j^- \le m$
\item $j\uparrow m$ if $j\not\le m$ and $j \le m^+$
\item $j \updownarrow m$ if $j \downarrow m$ and $j \uparrow m$.
\end{itemize}
\end{definition}

\begin{lemme} \mylabel{lem2-fleches} 
Let $L$ be a lattice. Then:
\begin{enumerate}
\item{$ \forall\ m \in M,\ \forall x \in L,\ (m \not\leq x\ 
\Longrightarrow\ \exists\ m \in M\ :\ x \leq m $ and $ j \uparrow m
$). }
\item{$\forall m \in M, \ \forall x \in L, \ (x \not\leq m
\Longrightarrow \exists j \in J: \  j \leq x$ and $j \downarrow m$).
  \mylabel{lem2-flechescond2}}
\end{enumerate}
Moreover, if $L$ is ULD, then Point \ref{lem2-flechescond2} becomes:\\
$$\forall\ m \in M,\ \forall x \in L,\ (x \not\leq m\ 
\Longrightarrow\ \exists\ j \in J\ :\ j \leq x \mbox{ and }  j \updownarrow m)$$
\end{lemme}

In a ULD lattice, for each $m\in M$, there are several $j\in J$ such that 
$j\updownarrow m$, but for each $j\in J$, there is exactly one $m\in M$ 
such that $j\updownarrow m$.
So $J$ is partitioned into $|M|$ sets 
$m_{\updownarrow}=\{j\in J, j\updownarrow m\}$.

Finally we give Birkhoff's representation theorem for distributive lattices:
\begin{theoreme}[Birkhoff] \cite{Bir33} \mylabel{thmbirkhoff}
A lattice is distributive if and only if it is isomorphic to the 
lattice of the ideals of the order induced on its meet-irreductibles.
\end{theoreme}

\section{Simple Chip Firing Games}
\mylabel{seccfgsimple}

\begin{figure}
\begin{center}
\begin{minipage}{5cm}
\begin{center}
\begin{picture}(30,40)
\node(a)(0,0){1}
\node(b)(30,0){1}
\node(c)(20,20){1}
\node(d)(5,30){0}
\gasset{ExtNL=y,NLdist=1}
\nodelabel[NLangle=180](a){$a$}
\nodelabel[NLangle=0](b){$b$}
\nodelabel[NLangle=0](c){$c$}
\nodelabel[NLangle=180](d){$d$}
\drawedge(a,c){}
\drawedge(b,c){}
\drawedge[curvedepth=2](c,d){}
\drawedge[curvedepth=-2](c,d){}
\end{picture}
\end{center}
\end{minipage}
\hspace{1cm}
\begin{minipage}{5cm}
\begin{center}
\includegraphics[scale=0.7]{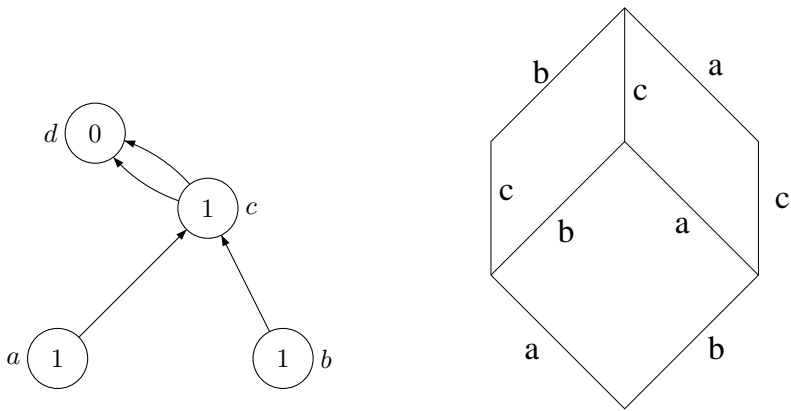}
\end{center}
\end{minipage}
\caption{A CFG and the lattice of its configuration space}
\label{figcfgtreillis}
\end{center}
\end{figure}

We will represent the lattice of the configuration space of a CFG by its Hasse diagram, and, when it is useful, we will label each edge between two configurations 
$\sigma$ and $\sigma'$ with the vertex $v$ such that $\sigma\stackrel{v}{\longrightarrow} \sigma'$.
Figure \ref{figcfgtreillis} shows an example of this representation.

\begin{note}
The order relation we have defined in the configuration spaces of CFGs is the dual of the one that is used in
 \cite{LP00}
(\emph{i.e.} in \cite{LP00} the initial configuration is the greatest element of the lattice, while it is the smallest element in this paper).
\cite{LP00} follows
 a convention for representing the states of dynamical models, that puts the initial configuration on top of the diagram.
We have chosen to do the opposite, following in that \cite{BLS91},
because this seems more natural in the context of order theory.
\end{note}

We give in this section a theorem that states that any CFG is equivalent to a \emph{simple} CFG, that is, a CFG where each vertex is fired at most once during an execution.
This result provides us with efficient tools for describing the configuration space of a CFG, thus greatly simplifiying the notations and proofs.
Moreover, new results can be derived quite simply from it.

\begin{definition}
A CFG is \emph{simple} if, during an execution, each vertex is fired at
most once.
\end{definition}

\begin{theoreme} \mylabel{thmsimple}
Any CFG that reaches a fixed point is equivalent to a simple CFG.
\end{theoreme}
Before giving the proof of this theorem, we first give the idea of the proof:
if a CFG is not simple, then
there exists at least one vertex $a$ that is
fired twice or more during an execution of $C$. We associate to $C$ a
CFG $C'$ equivalent to $C$, in which the vertex $a$ is split into two vertices
$a_0$ and $a_1$ that will be fired alternatively (the first firing of $a$ in $C$ corresponds to a firing of $a_0$ in $C'$, the second to a firing of $a_1$, 
and the $i$-th firing of $a$ corresponds to a firing of $a_0$ if $i$ is even, and to a firing of $a_1$ if $i$ is odd), 
so that each of them is fired less often than $a$ in $C$.
We acheive this by placing a large number of edges 
from $a_0$ to $a_1$ and from $a_1$ to $a_0$, and a large number of chips
in $a_0$, so that whatever configuration $C'$ is in, $a_1$ cannot 
contain enough chips to be fired
before $a_0$ is fired, thus bringing enough chips in $a_1$.
This large amount of chips will then move forth and back between the two vertices, guaranteeing that they will be fired alternatively.

\begin{proof}
Let $C$ be a non simple CFG with support graph $G=(V,E)$ and initial
configuration $\sigma$, and let $a$ be a vertex that is fired twice or more during an execution of $C$.
For a vertex $v$, we denote by $l(v)$ the number of loops on $v$.
We denote by $d{^>}_G(v)$ the number of edges going out of $v$ that are not loops (\emph{i.e.} $d{^>}_G(v)=d{_G}^+(v)-l(v)$).
We define dually $d{^<}_G(v)$.
The CFG $C'$ with support graph $G'=(V',E')$ and initial configuration
$\sigma'$ is defined in the following way:  let $N$ be twice the number of chips in $C$.
Then:
\begin{itemize}
\item let $V'=V\backslash \{a\} \cup \{a_0,a_1\}$, with $a_0\not\in V$ and
  $a_1\not\in V$.
\item $E'$ is defined by:
\begin{itemize}
    \item for each $v,w \in V\backslash \{a\}$, if there are $n$ edges$(v,w)$ in
      $E$, then there are $2n$ edges $(v,w)$ in $E'$.
    \item for each edge $(v,a)$ ($v\not= a$) in $E$, there is one
      edge $(v,a_0)$ 
      and one edge $(v,a_1)$ dans $E'$
    \item for each edge $(a,v)$ ($v\not= a$) in $E$, there are two
      edges $(a_0,v)$ 
      and two edges $(a_1,v)$ in $E'$
    \item for each loop $(a,a)$ in $E$, there is one loop
      $(a_0,a_0)$ and one loop $(a_1,a_1)$ in $E'$
    \item there are $N-d{_G}^>(a)$ edges both from $a_0$ to $a_1$ and
       from
      $a_1$ to $a_0$.
  \end{itemize}
\item for all $v\not=a$, $\sigma'(v)=2\sigma(v)$.\\
    $\sigma'(a_0)=\sigma(a)+N$, and $\sigma'(a_1)=\sigma(a)$.
\end{itemize}

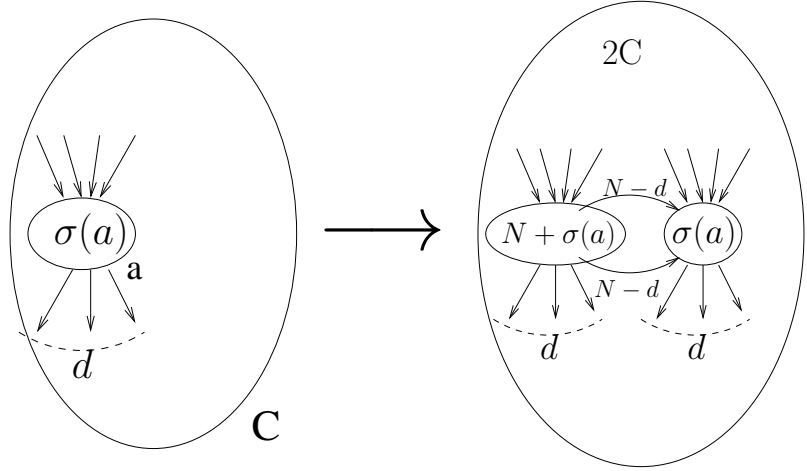
\begin{figure}
\begin{center}
\begin{minipage}{4cm}
\scalebox{0.75}{\input{figures/CFGnonsimpledetail.pstex_t}}
\end{minipage}
\begin{minipage}{2cm}
\scalebox{3}{$\longrightarrow$}
\end{minipage}
\begin{minipage}{4cm}
\scalebox{0.65}{\input{figures/CFGplussimple.pstex_t}}
\end{minipage}
\caption{Simplification of a CFG}
\label{figcfgsimple}
\end{center}
\end{figure}
Figure \ref{figcfgsimple} illustrates the construction.
We will prove  the following property: during an execution
of $C'$, every configuration of the game is such that one of the two
vertices $a_0$ or $a_1$ contains exactly $N$ chips more than the
other. This is true for the initial configuration. Since for each
$v\not= a_0,a_1$, there is the same number of edges from $v$ to $a_0$
than from $v$ to $a_1$, the firing of any other rule that one of the
$a_i$ does not change this property. Let us suppose now that we can
fire one of the rules $a_i$, for instance $a_0$
($a_0$ and $a_1$ can never be fired at the same time, because by construction there are not enough chips in the game).
 Let $N+x$ be the
number of chips in $a_0$ (there is then $x$ chips in $a_1$). The
outdegree of $a_0$ is $2d{_G}^>(a)+N-d{_G}^>(a)+l(a)=d{_G}^+(a)+N$. After
the firing of $a_0$, there are $x-d{_G}^>(a)$ chips in $a_0$, and
$N+x-d{_G}^>(a)$ in $a_1$. The property is thus verified.

We will now show that \confspace{C'} is isomorphic to \confspace{C}.
We start from the fact that the configuration space of the CFG obtained
from $C$ by doubling the initial configuration and the number of edges
of the support graph is isomorphic to the configuration space of
$C$. We will denote this CFG by $2C$.

The sum of the indegrees of $a_0$ and $a_1$ is equal to $2d{_G}^-(a)$,
i.e. the indegree of $a$ in $2C$. We can fire one of the rules $a_i$
if and only if it contains more than $N+d{_G}^+(a)$, that is if and
only if there is more than $N+2d{_G}^+(a)$ in the two vertices $a_0$
and $a_1$, which is $N$ chips more than the number of chips needed to fire $a$ in
$2C$. This firing will then give the same number of chips to the rest
of the graph as a firing of $a$ in $2C$. The other part of $C'$,
that is the rules different from $a_0$ or $a_1$, is as in $2C$ (except
the indegrees of the vertices $v$ such that there is an edge $(a,v)$
in $G$, but we have seen that this does not change the flow of chips
through the vertex). So the configuration space of $C'$ is isomorphic
to the configuration space of $2C$.

By this method we obtain a CFG $C'$ where the rules $a_0$ and $a_1$ are each
fired less often than in the initial CFG. By iterating this procedure,
we eventually obtain a simple CFG equivalent to the initial CFG.
\end{proof} 
We will now only consider simple CFGs, and we will also assume, without loss of generality, 
that their support graph has one and only one sink (denoted by $\bot$), so that the number of rules fired during an execution is equal to $|V|-1$ (therefore $|V|-1$ is also the height of the lattice of the configuration space).

The following result is due to \cite{LP00}:
\begin{lemme}
In a CFG, if, starting from the same configuration, two sequences of firing 
lead to the same configuration, then the set of rules fired in each sequence are the same.
\end{lemme}
This allows us to define the \emph{\shot{}} \sh{\sigma} of a configuration $\sigma$ as the set of rules fired to reach $\sigma$ from the initial configuration.
A subset $X\subseteq V\backslash\{ \bot\}$ is a \emph{valid} \shot{} if there exists a configuration $\sigma$ reachable from the initial configuration
such that $\sh{\sigma}=X$.
The lattice of the configuration space of a CFG is isomorphic to the lattice of the \shot{}s of its configurations ordered by inclusion \cite{LP00}.
The join is given by the following formula:
\begin{proposition} \mylabel{propjoincfg}
Let $C$ be a CFG, \confspace{C} its configuration space and $a,b$ two configurations.
The join of $a$ and $b$ in \confspace{C} is determined by:
$$\sh{a\vee b}=\sh{a}\cup\sh{b}$$
\end{proposition}

The following result appears with a different proof in \cite{BLS91}:
\begin{theoreme} \mylabel{thmcfglld}
The lattice of the configuration space of a CFG is ULD.
\end{theoreme}
\begin{proof}
Let $C$ be a CFG and $L=\lc$.
We will show that the interval between any element and the meet of its upper covers is a hypercube.
Let $x\in L$ and let $x_1,\ldots,x_n$ be its upper covers.
Each $x_i$ is obtained from $x$ by firing the vertex denoted by $i$, therefore $\sh{x_i}=\sh{x}\cup \{i\}$.
Now, the firing of a vertex $i$ does not prevent the firing of any vertex $i', i'\not=i$.
Indeed, when $i$ is fired, the number of chips in $i'$ stays the same or increases.
Therefore, for any subset $X$ of $\{1,\ldots,n\}$, the set $\sh{x}\cup X$ is a valid \shot{}, and so the interval between $x$ and $\bigvee\{x_1,\ldots,x_n\}$ is isomorphic to the set of all subsets of $\{1,\ldots,n\}$ and is a hypercube.
So by definition, $L$ is an ULD lattice.
\end{proof}

We will now show the link between the formulae characterising the elements 
and the join in a lattice of \lcfg{} (involving the \shot{}) and in a ULD lattice (involving the sets \M{x} as seen in Section \ref{secrecalls}).
Let $C$ be a CFG with support graph $G=(V,E)$ and $L=\lc$ be the lattice of its configuration space.
$L$ is ULD, so the height of $L$ is $|M|$.
We have seen that the height of $L$ is also equal to $|V|-1$.
We can talk equivalently of the rules of the CFG (the vertices that are fired during an execution, \ie all the vertices except the sink)
and of the meet-irreductibles of $L$ because
there is a bijection between $M$ and $V\backslash\{\bot\}$,
given by: $m \mapsto v$ if there exists an edge $(x,y)$ in $L$ such that $M_y\backslash M_x = m$ and $v$ is the rule that is applied to reach $y$ from $x$.
Moreover, this bijection preserves the formula for the join in the following way:
as seen in Section \ref{secrecalls}, the formula for the join in a ULD lattice is: $M_{x\vee y}=M_x\cap M_y$.
This is equivalent to: $M\backslash M_{x\vee y}=(M\backslash M_x)\cup(M\backslash M_y)$,
which is similar to the formula for the join given for the configuration space of a CFG in Proposition \ref{propjoincfg}: $\sh{C\vee C'}= \sh{C}\cup\sh{C'}$.
Therefore the \shot{} of a configuration $\sigma$ can be defined as $M\backslash M_{\sigma}$.

As an immediate consequence of Theorem \ref{thmsimple} we have:

\begin{corollaire} \mylabel{corintcfg}
Let $L$ be a lattice of \lcfg. Then every
interval of $L$ is also in \lcfg.
\end{corollaire}
\begin{proof}
Let us first recall that any interval of a lattice is a lattice.
Now let $L$ be a lattice of \lcfg{} and
$C$ be a corresponding CFG. 
The claim is true for any interval $[a,1_L]$ of $L$: 
let $\sigma$ be the configuration of $C$ that corresponds to $a$, then a
CFG with the same support graph as $C$ and with initial
configuration $\sigma$ has the lattice $[a,1_L]$ as configuration
space.

This result is also true for any interval $[0_L,b]$ in $L$. Indeed,
since $L$ is isomorphic to the configuration space of a simple CFG
$C'$, the element $b$ of $L$ partitions the vertices of $C'$ in two sets:
\sh{b} and $V\backslash \sh{b}$.
The interval $[0_L,b]$ is then the
configuration space of the CFG obtained from $C'$ by removing all the
edges going out of the vertices  in  $V\backslash \sh{b}$ (so that they
cannot be fired at all), and with the same initial configuration as
$C'$.

To conclude, simply notice that the interval $[a,b]$ in $L$ is the
intersection of the intervals $[0_L,b]$ et $[a,1_L]$.
\end{proof}

With the tools we have introduced in this section, we can describe efficiently a lattice in \lcfg{}, either (with the \shot{}s) by considering the CFG it is the configurations space of,
or (with the sets \M{x}) by use of the lattice theory.

\section{Study of the class \lcfg}
\mylabel{secstudy}
In this section we attempt to deepen the study of the class \lcfg.
We will show that it contains the class of distributive lattices, and that it is strictly included in the class of ULD lattices.

\begin{theoreme} \mylabel{thmdistrcfg}
Any distributive lattice is in \lcfg{}.
\end{theoreme}
\begin{proof}
Let $L$ be a distributive lattice and let $G=(M_L,<_L)$ be the graph
of the covering relation of the order induced on the meet-irreductibles of $L$.
Let $C$ be the CFG with support graph $G_C=(V,E)$ and initial configuration
$\sigma$, defined in the following way:
$V=M_L\cup \{\bot\}$, where $\bot\not\in M_L$, and
for each vertex $v$, the outgoing edges of $v$ in $E$ are  the edges going out from $v$ in $G$, plus:
\begin{itemize}
\item $d{_G}^-(v) - d{_G}^+(v)$  
   edges from $v$ to $\bot$ if $d{_G}^-(v) > d{_G}^+(v)$
\item  one edge from $v$ to $\bot$ if $v$ is isolated in $G$.
\end{itemize}
We notice that, for each vertex $v \in V$, $d_{G_C}^-(v) \le d_{G_C}^+(v)$.\\
The initial configuration is: for all $v\in V, \sigma(v)=d_{G_C}^+(v) -
 d_{G_C}^-(v)$.

$C$ verifies the following properties:
\begin{itemize}
\item each vertex can be fired once and only once (except the sink $\bot$ which
 is never fired).
Indeed, we can show by induction that the total number of chips that go out of
a vertex $v$ during an execution on $C$ is exactly $d^+(v)$: this is true for
each vertex without predecessor. Let now $v$ be a vertex such that this
property is verified for all of its predecessors. 
The immediate predecessors of $v$ can be fired exactly once, bringing 
$d_{G_C}^-(v)$ chips in $v$. By construction, 
$\sigma(v)=d_{G_C}^+(v) - d_{G_C}^-(v)$, so the property is verified for $v$.
\item A vertex $v$ can not be fired before all its predecessors
  in $G_C$ have been fired (notice that, since $G_C$ is constructed
  from the graph of an order, $G$ contains no cycle, therefore we do not reach a contradiction.
   Notice also that the set of predecessors of $v$ is the same 
   in $G$ and in $G_C$).
\end{itemize}
Therefore the \shot{} of each configuration  $\sigma$ is closed by the predecessor relation, and, since it does not contain $\bot$, it is an ideal of $M_L$.
Conversely, it is obvious that any ideal of $M_L$ is a valid \shot{} of $C$,
so the lattice of the configurations of $C$ is isomorphic to the lattice of the ideals of $M_L$, which by Birkhoff's Theorem (\ref{thmbirkhoff}) is isomorphic to $L$.
\end{proof}

\begin{theoreme}
Not all ULD lattices are the configuration space of a CFG.
\end{theoreme}
\begin{figure}
\begin{center}
\includegraphics[scale=0.5]{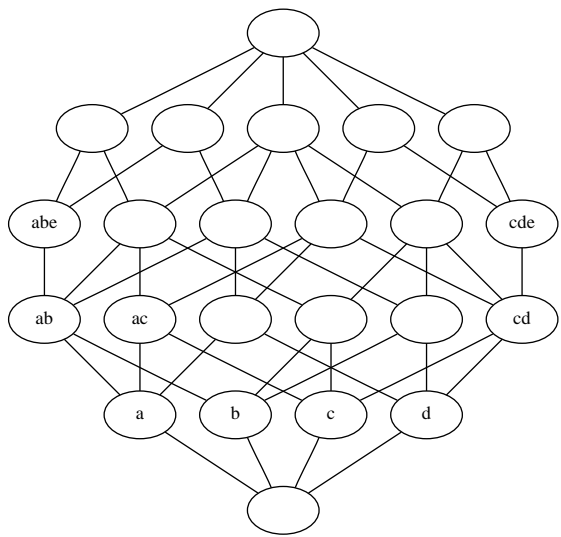}
\caption{A ULD lattice that is not the configuration space of a CFG}
\mylabel{figlldnoncfg}
\end{center}
\end{figure}
\begin{proof}
The lattice $L$  of Figure \ref{figlldnoncfg} is not in \lcfg.
Let us suppose that $L$ is in \lcfg{}.
By Theorem~\ref{thmsimple} it is the configuration space of a simple CFG $C$.
Since the height of $L$ is five, there are five distinct rules fired during an execution of $C$. 
At the beginning of the execution, four different rules can be fired: 
$a,b,c$ and $d$. The fifth rule, $e$, can be fired after the firing of either $a$ and $b$ or $c$ and $d$.
\begin{figure}
\begin{center}
\includegraphics[scale=0.5]{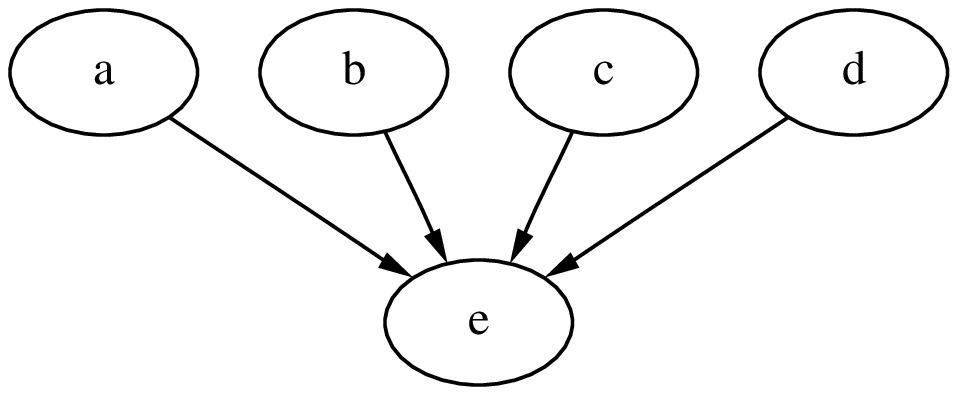}
\caption{}
\mylabel{figcfgpb}
\end{center}
\end{figure} 
The graph induced on $a,b,c,d$ and $e$ from the support graph of $C$ is
 therefore of the form displayed in Figure \ref{figcfgpb}.
Let $n$ be the number of chips that needs to fall into $e$ in order that $e$ can be fired. Let $n_a, n_b,n_c, n_d$ be the numbers of 
edges from $a,b,c$ and $d$ to $e$, that is the number of chips that fall into $e$ after the firing of one of these vertices.
We then have  $n_a+n_b\ge n$ and $n_c+n_d\ge n$. We can suppose that $n_a\ge n_b$ and $n_c\ge n_d$. 
But then $n_a+n_c\ge n$, and $e$ can be fired after the firing of $a$ and $c$.
So the configuration the \shot{} of which is $\{a,c\}$ should have three immediate successors:
$b$ and $d$, that can be fired since the beginning of the execution, and $e$.
Therefore we obtain a contradiction.
\end{proof}

\section{The coloured Chip Firing Game}
\mylabel{seccoloured}
We have seen in the previous section that the class \lcfg{} is strictly included in the class of ULD lattices.
We present now an extension of the Chip Firing Game, that generates exactly the class of ULD lattices.

For a graph $G=(V,E)$ and a set $X$ of colours, we call a \emph{coloured graph}
 the tuple $(V,E,X,col)$ where $col$ is a mapping from $E$ to $X$.
The \emph{restriction} of the graph to a colour $c\in X$ is the graph $(V,col^{-1}(c))$.
A \emph{coloured} CFG is defined over
a directed coloured multigraph $G=(V,E,X,col)$.
A configuration is given by a function $\sigma: V\rightarrow \N^{|X|}$
which associates to a vertex a number of chips of each colour.
Given a vertex $v$ and a colour $c$, we will denote by $\sigma_c(v)$ the number of chips of colour $c$ stored in $v$.
To each vertex is also associated
 a state function: at any time, a vertex can be 
 \emph{open}
 or \emph{closed}.
The evolution rule for this model is to \emph{open} a vertex.
One can open vertex $v$ if:
\begin{itemize}
\item $v$ is closed
\item
there exists a colour $c\in X$ such that $v$ can be fired (in the classical sense) in the restriction of the game to $c$
(that is, there are at least as many chips of colour $c$ in $v$ as there are edges of colour $c$ going out from $v$).
\end{itemize}
Opening a vertex consists in:
\begin{itemize}
\item marking it as open
\item
 for each colour $c$ in $X$,
 consider the restriction of the game to $c$ and to the set of open vertices, and play the game until the final configuration is reached.
\end{itemize}
 Notice that we have to ensure that the movements of chips that occur when opening a vertex stops after some time.
So we will consider
only graphs in which, for each colour $c$, 
the restriction of the game to $c$ is a (classical) convergent CFG
(this can be achieved by forbidding  
closed strongly connected component in the restriction of the graph to $c$ \cite{LP00}).

At the beginning of an execution, all the vertices are closed.
Since only closed vertices can be opened, coloured CFGs are convergent:
after some time, no vertex can be opened and the final configuration is reached.
They are also simple, therefore the configurations of the game are given by their \shot{}s.
The configuration space is ordered by the following relation:
$\sigma\le \sigma'\iff \sh{\sigma}\subseteq \sh{\sigma'}$.

The \emph{restriction} of a coloured CFG to a colour $c\in X$ is the game defined over the restriction of the support graph to $c$ 
such that, for each $v\in V$, the initial configuration is $\sigma_c(v)$.
The restriction of the game to a set of vertices is the game played on the induced subgraph with the corresponding restriction of the initial configuration.
In our figures, we will draw open vertices in gray.
The colours will be represented by numbers, the colour of an edge being indicated by its label.
In a vertex, a number $N_{c_1,\ldots,c_k}$ means that there are $N$ chips of colour $c_1$, $N$ chips of colour $c_k$, and so on, in the vertex.
For an example of execution of a coloured CFG see Figure \ref{figexecolore}

The coloured Chip Firing Game is an extension of the classical Chip Firing Game model.
Indeed we have the following result.
\begin{theoreme}
All convergent classical CFGs are equivalent to coloured CFGs
\end{theoreme}
\begin{proof}
By Theorem  \ref{thmsimple}, we know that a convergent CFG is equivalent to a simple CFG.
Now, any simple (convergent) CFG  can be viewed as a coloured CFG with only one colour.
\end{proof}

We will now show that the coloured CFGs generate exactly the class of ULD lattices.
The proof is given in two steps:
first we show in Theorem \ref{thmcoloredsuld} that it is included in the class of ULD lattices,
then we show in Theorem \ref{thmulddscolore} that it contains it.
We will use an intermediate theorem on lattice theory (Theorem \ref{thmidealquotient}) to prove Theorem \ref{thmulddscolore}.

\begin{theoreme} \mylabel{thmcoloredsuld}
The configuration space of a coloured CFG is a lattice.
Moreover, this lattice is ULD.
\end{theoreme}

\begin{proof}
Let $C$ be a coloured CFG. We will show that the set of the \shot{}s of the  reachable configurations of $C$ is closed under union. 
Let $\sigma_a$ and $\sigma_b$ be two reachable configurations of $C$, and let
$\sigma_o$ be a maximal configuration such that $\sigma_a$ and $\sigma_b$
can be reached from $\sigma_o$.
Let $A=\sh{\sigma_a},B=\sh{\sigma_b}$, and $O=\sh{\sigma_o}$.
 We clearly have $O\subseteq A$ and $O\subseteq B$. 
Let $a_1,\ldots,a_n \in A$ and $b_1,\ldots,b_n \in B$ be the vertices of $A$ and $B$ that can be opened in the configuration $\sigma_o$.
We have $\{a_1,\ldots,a_n\} \cap \{b_1,\ldots,b_n\} = \emptyset$ (otherwise, there exists a vertex $c\in A\cap B$ that can be opened in $\sigma_0$, 
leading to a configuration $\sigma_c$,
and both $\sigma_a$ and $\sigma_b$ can be reached from $\sigma_c$, 
which is greater than  $\sigma_o$).
To reach $\sigma_a$ from $\sigma_o$, we have to open all the vertices of $A\backslash O$. 
Doing this does not change the fact that the vertices $\{b_1,\ldots,b_n\}$ can be opened (the number of chips of each colour they contain stays the same or increases).
So the configuration $\sigma_{a'}$ such that $\sh{\sigma_{a'}} = A\cup \{b_1,\ldots,b_n\}$ is reachable from $\sigma_o$ and so from the initial configuration.
 Moreover, if $A'=\sh{\sigma_{a'}}$, then
$|(A\cup B) \backslash A'| = |(A\cup B) \backslash A| - n$.
By iterating this process, we eventually reach a configuration with \shot{} $A\cup B$.

The set of the \shot{}s of the reachable configuration has a smallest element, the emptyset, it is closed under union, so it is a lattice \cite{DP90}.
We will now show that this lattice is ULD:
let $\sigma_x$ be a configuration of $C$, $X=\sh{\sigma_x}$ and let $x_1,\ldots,x_n$ the vertices that can be opened in $\sigma_x$.
 Clearly, for any subset $Y$ of $\{\liste{x}\}$, the set $X\cup Y$ is a valid \shot{} of $C$ and the interval between $X$ and $X\cup \{\liste{x}\}$ is a hypercube, so the lattice is ULD.
\end{proof}
Given a coloured CFG $C$, we also denote by \confspace{C} the lattice of its configuration space.
Before giving the main theorem of this section, we need the following definitions:

\begin{definition}[$\sim$ relation]
Let $L$ be an ULD lattice.
We recall that $J$ is partitioned into $|M|$ sets 
$m_{\updownarrow}=\{j\in J, j\updownarrow m\}$.
We define $\sim$ by: $j\sim j'$ if and only if $j$ and $j'$ are in the same set $m_{\updownarrow}$.
\end{definition}

The $\sim$ relation induces an equivalence relation on the subsets of $J$:
two subsets $X$ and $Y$ are equivalent if and only if any element of $X\backslash Y$ is equivalent to an element of $Y$ and conversely.
We will consider the restriction of this relation on the ideals of $J$.
For each equivalence class, we define its representative element to be the maximal ideal in the class 
(it is unique because the union of ideals is an ideal, and if there were two maximal ideals in a class their union would still be in the class). 

\begin{figure}
\begin{center}
\includegraphics[scale=0.6]{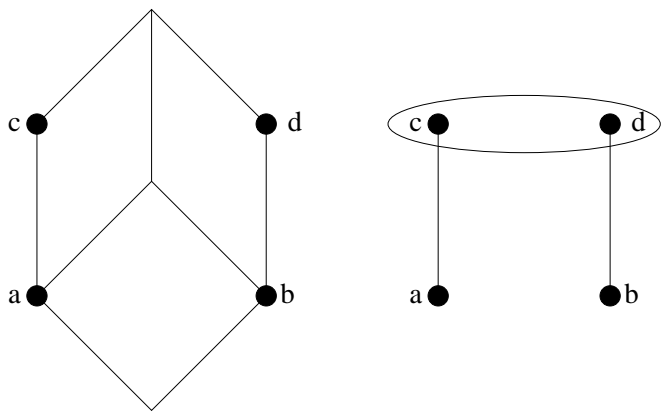}
\caption{An ULD lattice and the order on its join-irreductibles together with the $\sim$ relation}
\mylabel{figexordrejquot}
\end{center}
\end{figure}
Figure \ref{figexordrejquot} shows an exemple of the $\sim$ relation.
We have circled the join-irreductibles that are in the same equivalence class.

We then have the following result:
\begin{theoreme} \mylabel{thmidealquotient}
Let $L$ be an ULD lattice.
The set of the ideals on $J$ quotiented by $\sim$ is isomorphic to $L$.
\end{theoreme}
\begin{proof} 
We recall that $L$ is isomorphic to the set $\{J_x, x\in L\}$ ordered by inclusion.
We will proceed to show that, for any class of $\cal{O}(J)/\sim$, there exists $x\in L$ such that its representative element is equal to $J_x$,
then we will show that for any $x\in L$, $J_x$ is the representative element of some class of $\cal{O}(J)/\sim$.

In the sequel, we denote the set $M\backslash M_x$ by \sh{x} (following the notation \sh{x} introduced in Section \ref{seccfgsimple}).
We group here a few results that will be needed in the rest of the proof.
$m\in \sh{x}$ is equivalent by definition of \sh{x} to $x\not\le m$.
If $j\in J$ and $m\in M$ are such that $j\updownarrow m$, then
by definition of the $\updownarrow$ relation, we have $j\not\le m$, so $m\in\sh{j}$, and $j^-\le m$, so $m\not\in\sh{j^-}$.
Since $j\succ j^-$ we know that $|\M{j}\backslash\M{j^-}|=1$, so
we also have $\sh{j}=\sh{j^-}\sqcup\{m\}$ (where $\sqcup$ denotes the disjoint union).

Let $I\subseteq J$ be the representative element of a class of $\cal{O}(J)/\sim$, and let $x=\bigvee I$.
We will show that $I=J_x$. 
$I$ is obviously a subset of \J{x}, so $I\not= J_x$ implies that there exists $j\in J, j\in J_x$ (\ie $j\le x$ so $\sh{j}\subseteq\sh{x}$) and $j\not\in I$.
Let us assume that this is the case and
let $m$ be the meet-irreductible such that $j\updownarrow m$.
$m\in \sh{j}$ and $\sh{j}\subseteq\sh{x}$, so $m\in \sh{x}$.
We will prove that there exists $j'\in I, j\sim j'$.
Since $m\in \sh{x}$, and since $\sh{x}=\sh{\bigvee I}=\bigcup\{s(j), j\in I\}$, there exists $j''$ in $I$, $m\in \sh{j''}$,
 which means that $j''\not\le m$.
By Lemma \ref{lem2-fleches}, there exists $j'$ such that
$j'\updownarrow m$ and $j'\le j''$, so we have $j'\sim j$ and
 $j'\in I$ (because $j''\in I$ and $I$ is an ideal).

So any element of $J_x\backslash I$ is equivalent to some element of $I$, and since $I\subset J_x$,
 $I\sim J_x$.
So $I$ is not the maximal ideal in its class, which is a contradiction.
Therefore all the representative elements of $\cal{O}(J)/\sim$ are equal to a set $J_x$, for some $x\in L$.

We will now show that for all $x\in L$, there is a class of $\cal{O}(J)/\sim$ such that $J_x$ is its representative element.
Since (as we have just seen) the representative element of any class is equal to $J_y$, for some $y\in L$, the only way for $J_x$ not to be the representative element of a class is for $J_x$ to be included in a set $J_y$, with $J_x\sim J_y$.
We will prove that this is not possible, \emph{i.e.} that $J_x\subseteq J_y$ and $J_x\sim J_y$ implies  $J_x = J_y$.
Let us assume that there exists $y$ such that $J_x\subseteq J_y$ and $J_x\sim J_y$.
Let $j$ be a minimal element of $J_y\backslash J_x$, so $J_x\cup \{j\}$ is an ideal of $J$.
We will show that $j\le x$, obtaining thus a contradiction.
Since $J_x\sim J_y$, there exists $j'\in J_x, j\sim j'$.
Since $J_x\cup \{j\}$ is an ideal of $J$, all the strict predecessors of $j$ in $J$ are elements of \J{x}.
Since $j^-$ is the only immediate predecessor of $j$, any strict predecessor of $j$ in $L$ is a predecessor of $j^-$, so that $J_{j^-}=J_j\backslash \{j\}$.
Therefore $\J{j^-}\subseteq \J{x}$ (which means $j^-\le x$).
Since $j\sim j'$, there exists $m\in M$ such that $j\updownarrow m$ and $j'\updownarrow m$.
We have seen that this means that $m\in \sh{j'}$ and that $\sh{j}=\sh{j^-}\cup \{m\}$.
Therefore $\sh{j}\subseteq \sh{j^-}\cup \sh{j'}$, so $j\le j^-\vee j'$.
We already know that $j^-\le x$ and $j'\le x$, so $j\le x$, which is impossible because we have assumed that $j\not\in J_x$.

So the representative elements of the equivalence classes of $\cal{O}(J)/\sim$ are exactly the sets $J_x$, so $L$ is isomorphic to the set of the ideals of $J$ quotiented with respect to the $\sim$ relation.
\end{proof}

\begin{note}
This theorem is quite close to Nourine's work on coloured ideals \cite{nou00}.
In his work Nourine defines a coloured ideal to be the set of colours associated to an ideal of a coloured poset (\emph{i.e.} a poset where a colour $c(x)$ is associated to each element $x$ in such a way that $x<y$ implies $c(x)\not=c(y)$).
This implies the following results:
the set of the coloured ideals of a coloured poset is an ULD lattice,
and, given an ULD lattice $L$, there exists a coloured poset such that the lattice of its coloured ideals is $L$.
Theorem \ref{thmidealquotient} is similar to this last result where the coloured poset is defined from $J$, and where a different colour is associated to each equivalence class of the $\sim$ relation.
\end{note}

\begin{theoreme} \mylabel{thmulddscolore}
Let $L$ be an ULD lattice.
Then there exists a coloured CFG $C$ such that $L=\confspace{C}$.
\end{theoreme}
\begin{proof} 
Let $L$ be an ULD lattice.
We are going to build a coloured CFG
 $C$ with support graph $(V,E,X,col)$
satisfying $L=\confspace{C}$.
We will first construct a coloured CFG $\tilde{C}$ 
such that \confspace{\tilde{C}} is the lattice of the ideals of $J$.
The support graph of $\tilde{C}$ is $(J\cup\{\bot\},\tilde{E},X,\tilde{col})$,
and $\tilde{C}$ is defined in the following way:
we define $X$ such that $|X|=|J|$. To each element $j$ of $J$ we associate in a bijective way a colour $c_j$.
For each element $j$ of $J$ we define $\downarrow j$ to be the smallest ideal that contains $j$.
Then for each $j\in J$, we construct with the colour $c_j$ a (classical) CFG the configuration space of which is the lattice of the ideals of $\downarrow j$.
We do this by
applying the algorithm descibed in the proof of Theorem \ref{thmdistrcfg}.

\begin{figure}
\begin{center}
\begin{picture}(80,50)
\node(a)(30,0){$1_{1,3}$}
\node(b)(50,0){$1_{2,4}$}
\node(c1)(30,15){}
\node(c2)(50,15){}
\node(puits)(40,30){}
\drawedge(a,c1){$1$}
\drawedge(b,c2){$2$}
\drawedge(c1,puits){1}
\drawedge[ELside=r](c2,puits){$2$}
\drawqbedge(a,10,15,puits){$3$}
\drawqbedge[ELside=r](b,70,15,puits){$4$}
\end{picture}
\end{center}
\caption{The coloured CFG $\tilde{C}$ for the lattice of Figure \ref{figexordrejquot}}
\mylabel{figexcfgeclate}
\end{figure}
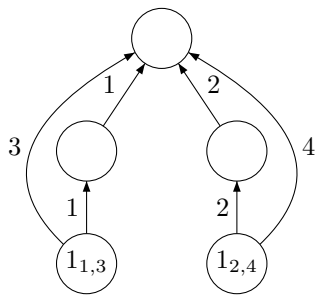
Figure \ref{figexcfgeclate} shows an example of the construction of $\tilde{C}$.

$C$ is obtained by contracting $\tilde{C}$ with respect to the $\sim$ relation:
$V=J/\sim\cup\{\bot\}$ (remark that $V$ is in bijection with $m\cup\{\bot\}$), and for each $v\in V$ and each colour $c\in X$, we have (if $j_1,\ldots,j_k$ are all the vertices in the class that $v$ represents):
\begin{itemize}
\item $succ_c(v)=\left( \bigcup succ_c(j_i)\right)/\sim$
\item $pred_c(v)=\left( \bigcup pred_c(j_i)\right)/\sim$
\end{itemize}
and the initial configuration is given by: $\sigma_c(v)=\sum\sigma_c(j_i)$.

We will now prove that \confspace{C} is isomorphic to the set of the ideals of $J$ quotiented by the $\sim$ relation.
To a \shot{} $s$ of $C$ we associate the ideal $I(s)$ of $J$ that is the greatest ideal included in $\bigcup\{m_{\updownarrow},m\in s\}$
(notice that $I(s)$ is the representative element of its class).

We will now show that if $s\not=s'$ are two valid \shot{}s of $C$, then $I(s)\not=I(s')$.
If $s\not=s'$ are two valid \shot{}s, then
without loss of generality, we can assume that there exists $m\in s, m\not\in s'$.
If $m$ has been opened, it is because it has gathered a sufficient number of chips of a given colour $c$. Let $j$ be the join-irreductible such that $c=c_j$.
By construction of $C$, there exists $j'\in m_{\updownarrow}, j'\le j$.
If $m$ has been opened for the colour $c$, then all the vertices corresponding to the predecessors of $j'$ must have been opened before,
so $\downarrow j'\subseteq \bigcup\{m_{\updownarrow},m\in s\}$,
which means that $j'\in I(s)$. Since it is obvious that $j'\not\in I(s')$, 
we have $I(s)\not=I(s')$.
By construction, if $s\subseteq s'$, $I(s)\subseteq I(s')$.

Conversely, if $I\in \cal{O}(J)/\sim$, we will show that there exists a valid \shot{} $s$ such that $I=I(s)$ (notice that there is a unique such $s$).
Let $M(I)=\{m\in M,\exists j\in J, j\in m_{\updownarrow}\}$.
We will show that $M(I)$ is a valid \shot{}.
Let $j_i,\ldots,j_k$ be a linear extension of $I$.
For each $i$, we define $m_i$ to be the class $j_i$ is in.
We will construct a sequence $s_1,\ldots,s_k$ of valid \shot{}s such that $s_i\subseteq s_{i+1}$, and, for each $i$, $M(\{j_1,\ldots,j|i\})=s_i$.
$s_1=m_1$ is a valid \shot{}, because $j_1$ is minimal in $J$, so $m_1$ can be opened for the colour $c_{j_1}$.
If $s_i=\bigcup_{l=1}^i [m_j]$ is a valid \shot{}, then $s_{i+1}=\bigcup_{l=1}^{i+1} [m_i]$, because all the predecessors of $j_{i+1}$ are open in $s_i$, and by construction of $C$, $m_{i+1}$ can be opened for the colour $c_{i+1}$ if all the predecessors of $j_{i+1}$ are open.
So $s_k=M(I)$, which means that $I(s_k)=I$.
In the same way, if $I\subseteq I'$ are two elements of $ \cal{O}(J)/\sim$, any linear extension of $I$ is the beginning of a linear extension of $I'$, so 
the \sh{} $s$ such that $I(s)=I$ is included in the \sh{} $s'$ such that $I(s')=I'$.
\end{proof}

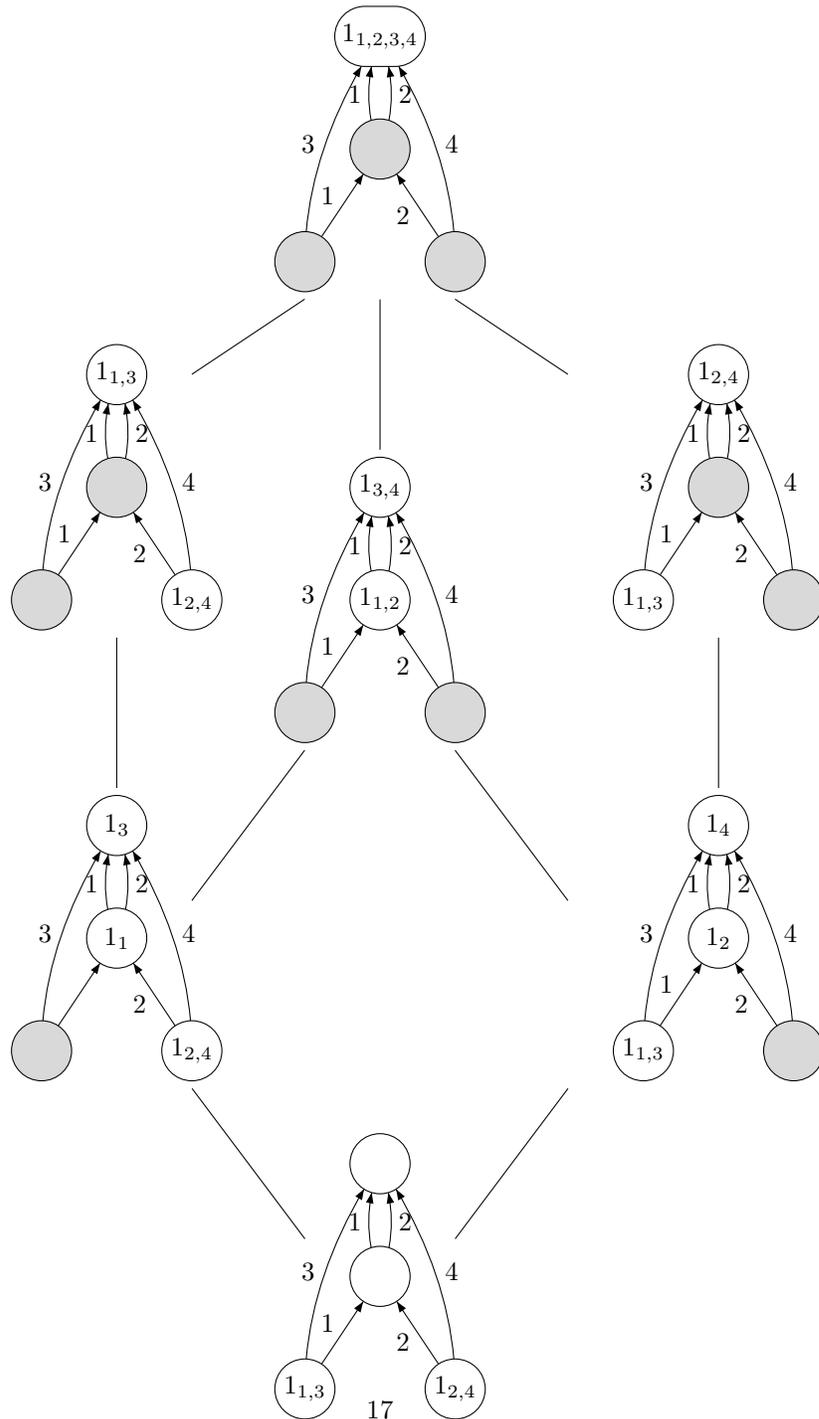
\begin{figure} 
\begin{center}
\begin{picture}(70,200)
\put(35,0){        
\node(a)(-10,0){$1_{1,3}$}
\node(b)(10,0){$1_{2,4}$}
\node(c)(0,15){}
\node(puits)(0,30){}
\drawedge(a,c){$1$}
\drawedge(b,c){$2$}
\drawqbedge(c,-3,22,puits){$1$}
\drawqbedge[ELside=r](c,3,22,puits){$2$}
\drawqbedge(a,-10,15,puits){$3$}
\drawqbedge[ELside=r](b,10,15,puits){$4$}}

\put(0,45){                
\node[fillgray=0.85](a)(-10,0){}
\node(b)(10,0){$1_{2,4}$}
\node(c)(0,15){$1_1$}
\node(puits)(0,30){$1_3$}
\drawedge(a,c){}
\drawedge(b,c){$2$}
\drawqbedge(c,-3,22,puits){1}
\drawqbedge[ELside=r](c,3,22,puits){$2$}
\drawqbedge(a,-10,15,puits){$3$}
\drawqbedge[ELside=r](b,10,15,puits){$4$}}

\put(80,45){              
\node(a)(-10,0){$1_{1,3}$}
\node[fillgray=0.85](b)(10,0){}
\node(c)(0,15){$1_2$}
\node(puits)(0,30){$1_4$}
\drawedge(a,c){$1$}
\drawedge(b,c){$2$}
\drawqbedge(c,-3,22,puits){1}
\drawqbedge[ELside=r](c,3,22,puits){$2$}
\drawqbedge(a,-10,15,puits){$3$}
\drawqbedge[ELside=r](b,10,15,puits){$4$}}

\put(35,90){                 
\node[fillgray=0.85](a)(-10,0){}
\node[fillgray=0.85](b)(10,0){}
\node(c)(0,15){$1_{1,2}$}
\node(puits)(0,30){$1_{3,4}$}
\drawedge(a,c){$1$}
\drawedge(b,c){$2$}
\drawqbedge(c,-3,22,puits){1}
\drawqbedge[ELside=r](c,3,22,puits){$2$}
\drawqbedge(a,-10,15,puits){$3$}
\drawqbedge[ELside=r](b,10,15,puits){$4$}}

\put(35,150){                  
\node[fillgray=0.85](a)(-10,0){}
\node[fillgray=0.85](b)(10,0){}
\node[fillgray=0.85](c)(0,15){}
\node[Nadjust=w](puits)(0,30){$1_{1,2,3,4}$}
\drawedge(a,c){$1$}
\drawedge(b,c){$2$}
\drawqbedge(c,-3,22,puits){1}
\drawqbedge[ELside=r](c,3,22,puits){$2$}
\drawqbedge(a,-10,15,puits){$3$}
\drawqbedge[ELside=r](b,10,15,puits){$4$}}

\put(80,105){                 
\node(a)(-10,0){$1_{1,3}$}
\node[fillgray=0.85](b)(10,0){}
\node[fillgray=0.85](c)(0,15){}
\node(puits)(0,30){$1_{2,4}$}
\drawedge(a,c){$1$}
\drawedge(b,c){$2$}
\drawqbedge(c,-3,22,puits){1}
\drawqbedge[ELside=r](c,3,22,puits){$2$}
\drawqbedge(a,-10,15,puits){$3$}
\drawqbedge[ELside=r](b,10,15,puits){$4$}}

\put(0,105){                
\node[fillgray=0.85](a)(-10,0){}
\node(b)(10,0){$1_{2,4}$}
\node[fillgray=0.85](c)(0,15){}
\node(puits)(0,30){$1_{1,3}$}
\drawedge(a,c){$1$}
\drawedge(b,c){$2$}
\drawqbedge(c,-3,22,puits){1}
\drawqbedge[ELside=r](c,3,22,puits){$2$}
\drawqbedge(a,-10,15,puits){$3$}
\drawqbedge[ELside=r](b,10,15,puits){$4$}}

\gasset{AHnb=0}
\drawline(25,20,10,40)
\drawline(45,20,60,40)
\drawline(60,65,45,85)
\drawline(10,65,25,85)
\drawline(35,125,35,145)
\drawline(80,80,80,100)
\drawline(0,80,0,100)

\drawline(60,135,45,145)
\drawline(10,135,25,145)
\end{picture}

\caption{The configuration space of the coloured CFG obtained from the CFG $\tilde{C}$ of Figure \ref{figexcfgeclate}}
\mylabel{figexecolore}
\end{center}
\end{figure}

Figure \ref{figexecolore} shows an example of the execution of a coloured CFG obtained as described in the proof.

\section*{Conclusion}
In this paper we have attempted to define exactly the class \lcfg{}
of lattices that are the configuration space of a CFG.
We have shown that this class is strictly between (with respect to inclusion) two well-known classes, the distributive and the ULD lattices.
We have also presented an extension of the model, the coloured CFG, which increases the number of lattices that can be represented.
It is a natural extension since any classical CFG can be seen as a coloured CFG.
The converse, that is, transform a coloured CFG (when it is possible) into a classical CFG, and the characterisation of the class \lcfg{} (with, for instance, an algorithm that, given a lattice, constructs a corresponding CFG or fails if the lattice is not in \lcfg{}), remains to be done.

\section*{Acknowlegments}
The authors wish to thank Matthieu Latapy for the help and the useful comments he gave for this paper.

\clearpage

\bibliographystyle{alpha}
\bibliography{../../biblio}
\end{document}

%% file: figures/CFGnonsimpledetail.pstex_t
\begin{picture}(0,0)%
\includegraphics{CFGnonsimpledetail.pstex}%
\end{picture}%
\setlength{\unitlength}{3947sp}%
\begingroup\makeatletter\ifx\SetFigFont\undefined%
\gdef\SetFigFont#1#2#3#4#5{%
  \reset@font\fontsize{#1}{#2pt}%
  \fontfamily{#3}\fontseries{#4}\fontshape{#5}%
  \selectfont}%
\fi\endgroup%
\begin{picture}(2416,3612)(593,-2767)
\put(1126,-2161){\makebox(0,0)[lb]{\smash{\SetFigFont{20}{24.0}{\rmdefault}{\mddefault}{\updefault}\special{ps: gsave 0 0 0 setrgbcolor}$d$\special{ps: grestore}}}}
\put(976,-1036){\makebox(0,0)[lb]{\smash{\SetFigFont{20}{24.0}{\rmdefault}{\mddefault}{\updefault}\special{ps: gsave 0 0 0 setrgbcolor}$\sigma(a)$\special{ps: grestore}}}}
\end{picture}

%% file: figures/CFGplussimple.pstex_t
\begin{picture}(0,0)%
\includegraphics{CFGplussimple.pstex}%
\end{picture}%
\setlength{\unitlength}{3947sp}%
\begingroup\makeatletter\ifx\SetFigFont\undefined%
\gdef\SetFigFont#1#2#3#4#5{%
  \reset@font\fontsize{#1}{#2pt}%
  \fontfamily{#3}\fontseries{#4}\fontshape{#5}%
  \selectfont}%
\fi\endgroup%
\begin{picture}(3166,4512)(818,-3817)
\put(2026, 89){\makebox(0,0)[lb]{\smash{\SetFigFont{20}{24.0}{\rmdefault}{\mddefault}{\updefault}\special{ps: gsave 0 0 0 setrgbcolor}2C\special{ps: grestore}}}}
\put(2851,-2761){\makebox(0,0)[lb]{\smash{\SetFigFont{20}{24.0}{\rmdefault}{\mddefault}{\updefault}\special{ps: gsave 0 0 0 setrgbcolor}$d$\special{ps: grestore}}}}
\put(1426,-2761){\makebox(0,0)[lb]{\smash{\SetFigFont{20}{24.0}{\rmdefault}{\mddefault}{\updefault}\special{ps: gsave 0 0 0 setrgbcolor}$d$\special{ps: grestore}}}}
\put(2026,-1186){\makebox(0,0)[lb]{\smash{\SetFigFont{14}{16.8}{\rmdefault}{\mddefault}{\updefault}\special{ps: gsave 0 0 0 setrgbcolor}$N-d$\special{ps: grestore}}}}
\put(1951,-2161){\makebox(0,0)[lb]{\smash{\SetFigFont{14}{16.8}{\rmdefault}{\mddefault}{\updefault}\special{ps: gsave 0 0 0 setrgbcolor}$N-d$\special{ps: grestore}}}}
\put(1051,-1636){\makebox(0,0)[lb]{\smash{\SetFigFont{17}{20.4}{\rmdefault}{\mddefault}{\updefault}\special{ps: gsave 0 0 0 setrgbcolor}$N+\sigma(a)$\special{ps: grestore}}}}
\put(2701,-1636){\makebox(0,0)[lb]{\smash{\SetFigFont{20}{24.0}{\rmdefault}{\mddefault}{\updefault}\special{ps: gsave 0 0 0 setrgbcolor}$\sigma(a)$\special{ps: grestore}}}}
\end{picture}